\documentclass[reqno, oneside]{amsart}
\usepackage{style}

\begin{document}

\title[Analytic approximation in $L^p$]
{Analytic approximation in $L^p$ and coinvariant subspaces of the Hardy space}
\author{R.~V.~Bessonov}
\address{St.Petersburg Department of Steklov Mathematical Institute RAS ({\normalfont 27, Fontanka, St.Petersburg, 191023, Russia}) and St.Petersburg State University ({\normalfont \hbox{7-9}, Universitetskaya nab., St.Petersburg, 199034, Russia})}
\email{bessonov@pdmi.ras.ru}
\thanks{This work is partially supported by the RFBR grants 12-01-31492, 11-01-00584, by Moebius Contest Foundation for Young Scientists and by the Chebyshev Laboratory  (Department of Mathematics and Mechanics, St. Petersburg State University) under RF Government grant 11.G34.31.0026}
\date{}

\begin{abstract}
We generalize a classical result by A.\ Macintyre and W.\ Rogosinski on best $H^p$--approximation in $L^p$ of rational functions. For each inner function $\theta$ we give a description of $H^p$--badly approximable functions in $\bar \theta H^p$.
\end{abstract}

\maketitle

\section{Introduction}
The classical problem of best analytic approximation in $L^p$ on the unit circle $\T$ reads as follows: given a function $g \in L^p$, find a function $p_g$ in the Hardy space $H^p$ such that
$$
\|g - p_g\|_{L^p} = \dist_{L^p} (g, H^p).
$$
In 1920, F.Riesz proved \cite{MR1555162} that best $H^1$--approximation in $L^1$ of a trigonometric polynomial of degree $n$ is an analytic polynomial of degree at most~$n$. His result was generalized in 1950 by A.\ Macintyre and W.\ Rogosinski \cite{MR0036314}, who treat the problem of best analytic approximation in $L^p$ for rational functions with finite number of poles in the open unit disk.  
\begin{Thm}[A.\ Macintyre, W.\ Rogosinski] \label{mrt}
Let $1 \le p \le \infty$, and let $g$ be a rational function with $n$ poles $\beta_i$ in $|z|<1$, each counted according to multiplicity. Then best $H^p$--approximation $p_g$ of the function $g$ exists uniquely. Moreover, there exist $n-1$ numbers $\alpha_i$ with $|\alpha_i| \le 1$ such that 
\begin{equation}\label{e1}
g - p_g = const \cdot \prod{}^{'} \frac{z - \alpha_i}{1 - \bar \alpha_i z} \prod_{1}^{n-1} (1 - \bar \alpha_i z)^{2/p}  \prod_{1}^{n} \frac{1 - \bar \beta_i z}{z - \beta_i} (1 - \bar \beta_i z)^{-2/p},  \\
\end{equation}
where $\prod{}^{'}$ is extended over all, some, or none of the $\alpha_i$ with $|\alpha_i| < 1$.
\end{Thm}
Among other things, this result shows that best 
$H^1$--approximation of a rational function is a rational function as well. The same holds for best $H^\infty$--approximation. 

\smallskip

In 1953, W.\ Rogosinski and H.\ Shapiro \cite{MR0059354} presented a uniform approach to the problem of best analytic approximation in $L^p$ based on duality for classes $H^p$. Their paper contains a refined (but still rather complicated) proof of Theorem \ref{mrt}.

\smallskip

The matrix-valued case of the problem of best analytic approximation has been studied extensively in the last years.
In particular, V.Peller and V.Vasyunin \cite{PV} consider this problem for rational matrix-valued functions motivated by applications in $H^\infty$--\! Control Theory. A survey of results related to best analytic approximation in $L^p$ of matrix-valued functions can be found in L.Baratchart, F.Nazarov, V.Peller~\cite{BNP}.

\smallskip

Our aim in this note is to give a short proof of Theorem \ref{mrt} and present its analogue in more general situation. We will consider the problem of best analytic approximation for functions of the form $h/\theta$, where $h \in H^p$ and $\theta$ is an inner function. If $\theta$ is a finite Blaschke product we are in the setting of Theorem \ref{mrt}. In general, functions of the form $h/\theta$ may have much more complex behaviour near the unit circle than rational functions. To be more specific, we need some definitions.

A bounded analytic function $\theta$ in the open unit disk is called inner if $|\theta| = 1$ almost everywhere on the unit circle $\T$ in the sense of angular boundary values. Given an inner function $\theta$, define the coinvariant subspace $\Kthp$ of the Hardy space $H^p$ by the formula $\Kthp = H^p \cap \bar z \theta  \ov{H^p}$. Here and in what follows we identify the Hardy space $H^p$ in the open unit disc $\D$ with the corresponding subspace of the space $L^p$ on the unit circle $\T$ via angular boundary values. All the information we need about Hardy spaces is available in Sections II and IV of \cite{MR628971}. Basic theory of coinvariant subspaces $\Kthp$ can be found in \cite{MR827223}, \cite{MR1289670}.

\medskip

Our main result is the following.

\begin{Thm}\label{t1}
Let $\theta$ be an inner function and let $1 \le p < \infty$. Take a function 
$g \in \bar \theta H^p$ and denote by $p_g$ its best $H^p$--approximation. The function $g - p_g$ can be uniquely represented in the form $g - p_g =  c \cdot  \bar \theta I F^{2/p}$, where $c = \dist_{L^p}(g , H^p)$, $F$ is an outer function in $\Kth$ of unit norm, $F(0)>0$, and $I$ is an inner function such that $I F \in \Kth$. 
\end{Thm}

Taking $\theta = z^n$ and $p=1$ in Theorem \ref{t1}, we get the mentioned result by F.Riesz: trigonometric polynomials are preserved under best analytic approximation in $L^1$. Our paper contains the fourth proof of this fact (see Section \ref{s3}); previous proofs can be found in \cite{MR1555162}, \cite{MR0036314}, \cite{MR0382963}. Similarly, Theorem \ref{mrt} follows from Theorem \ref{t1} by taking the function $\theta$ to be a finite Blaschke product. The choice $\theta=e^{iaz}$ leads to the following fact.

\begin{Thm}\label{t2}
Let $g$ be a function in $L^1(\R)$ with compact support of Fourier transform: $\supp \hat g \subset [-a, a].$ Then we have $\supp \hat p_g \subset [0, a]$ for best $H^1$--approximation~$p_g$ of the function $g$ .
\end{Thm}
Proofs of Theorems \ref{mrt}, \ref{t1}, \ref{t2} are given in Sections  \ref{s3}, \ref{s2}, \ref{s4}, correspondingly. In Section \ref{s5} we discuss how the problem of best analytic approximation in $L^p$ for functions from $\bar \theta H^p$ can be reduced to a special problem of interpolation. 

\section{Proof of Theorem \ref{t1}}\label{s2}
We need the following known result from \cite{MR1097956}.
\begin{Lem}[K.Dyakonov]\label{l1}
A nonnegative function $\phi$ can be represented in the form $\phi =|F|^2$ for some outer function $F \in \Kth$ if and only if $\phi \in z \bar \theta H^1$.
\end{Lem}

The proof is included for completeness.

\beginpf Let $\phi$ be a function of the form $\phi =|F|^2$, where $F \in H^2 \cap \bar z \theta \ov{H^2}$. Take a function $G \in H^2$ such that $F = \bar z \theta \bar G$. We have $\phi = z \bar \theta GF \in z \bar \theta H^1$, as required. 

Conversely, consider a nonnegative function $\phi \in z \bar \theta H^1$. Since $\theta$ is unimodular on the unit circle $\T$, we have $\log\phi \in L^1$. Let $F$ be the outer function in $H^2$ with modulus $\sqrt{\phi}$ on $\T$. We have $\bar z \theta |F|^2 \in H^1$. Hence, $\bar z \theta |F|^2 = IF^2$ for an inner function $I$. Thus, the function $F =  \bar z \theta \ov{IF}$ belongs to the subspace $\bar z \theta \ov{H^2}$. It follows that $F \in \Kth$, which completes the proof.  \qed

\bigskip

\noindent {\bf Proof of Theorem 2.} Let $g$ be a function in the subspace $\bar \theta H^p$, where $1 \le p <\infty$. Denote by $p'$ the conjugate exponent to $p$. There exist functions $p_g \in H^p$, $h_g \in zH^{p'}$ satisfying
\begin{equation}\label{e3}
\|g - p_g\|_{L^p} = \dist_{L^p}(g, H^p) = \int_{\T} (g-p_g)h_g\,dm, \quad \|h_g\|_{L^{p'}} = 1,
\end{equation}
where $m$ denotes the normalized Lebesgue measure on $\T$. This well-known fact was first established in \cite{MR0059354}; its modern proof can be found, e.g., in Section IV of~\cite{MR628971}. 
Denote by $f$ the function $g - p_g \in \bar \theta H^p$ and set $c = \|f\|_{L^p} = \dist_{L^p} (g, H^p)$. It follows from \eqref{e3} that we have equality in the H\"{o}lder inequality $\|fh_g\|_{L^1} \le \|f\|_{L^p}\|h_g\|_{L^{p'}}$. Therefore, $fh_g = c^{1-p}\cdot|f|^{p}$.

The function $fh_g$ belongs to the subspace $z \bar \theta H^1$. Hence, the function $|f|^{p}$ belongs to $z \bar \theta H^1$ as well, and we see from Lemma \ref{l1} that $|f|^{p} = c^{p} |F|^2$ for an outer function $F \in \Kth$ of unit norm. We may assume that $F(0) > 0$. The function $\theta f$ lies in $H^p$ and has modulus $c |F|^{2/p}$. It follows that $\theta f = c I F^{2/p}$ for an inner function $I$. Let us prove that $IF \in \Kth$. By the construction, we have 
\begin{equation}\label{e4}
c^p|F|^2 = |f|^p = c^{p-1} f h_g = c^p\cdot \bar \theta I F^{2/p} h_g.
\end{equation}
Hence, the function $h_g \in z H^{p'}$ has the form $h_g = z J F^{2/{p'}}$ where $J$ is an inner function. From \eqref{e4} we get the formula $z \bar \theta IJ F = \bar F$. This yields the fact that $IF \in \bar z \theta \ov{H^2}$. Thus, the inclusion $IF \in \Kth$ is proved. By the construction,  $f = c \cdot \bar \theta I F^{2/p}$.

Now prove that functions $I,F$ in the statement of the theorem are determined uniquely. For $1 \le p < \infty$, best $H^p$--approximation $p_g$ of the function $g$ is unique, see \cite{MR0059354} or Section IV in \cite{MR628971}. Hence, the function $c \cdot I F^{2/p} = \theta (g - p_g)$ is determined uniquely. It remains to use  uniqueness in the inner-outer factorization for functions in $H^p$.   
\qed

\medskip

\noindent {\bf Remark 2.1.} In the case $p= \infty$, Theorem \ref{t1} holds provided the dual extremal function $h_g \in z H^1$ in formula \eqref{e3} exists. Indeed, under this assumption best $H^\infty$--approximation $p_g$ is unique and we get from \eqref{e3} that $f h_g = c|h_g|$, where $f = g - p_g$ and $c = \|f\|_\infty = \dist_{L^\infty}(f, H^\infty)$. As above, there exists an outer function $F \in \Kth$  such that $|h_g| = |F|^2$. Hence, $f h_g = c|F|^2$ and we have $z \bar \theta IF^2 = c|F|^2$ for some inner function $I$. It follows that $IF \in \Kth$ and $f = c\bar \theta I$, as required. 

It can be shown that the dual extremal function $h_g$ exists for every continuous function $g$ on the unit circle $\T$, see \cite{MR0049322} or Section IV in \cite{MR628971}. In particular, it exists for every rational function with poles in the open unit disk. This will allow us to prove Theorem~\ref{mrt} in the case $p=\infty$, see details in the next section.

\medskip

\noindent {\bf Remark 2.2.} As we have seen in the proof of Theorem \ref{t1}, the dual extremal function $h_g$ to the function $g$ is given by the formula $h_g = z J F^{2/p'}$, where $J$ is the inner function such that 
$IJF = \bar z \theta \bar F$.  It can be shown that every inner function $U$ for which $UF \in \Kth$ is a divisor of the function $IJ$, see Theorem 2 in \cite{MR1097956}.

\section{Proof of Theorem \ref{mrt}}\label{s3}
Let us first prove the classical result  by F.Riesz on best analytic approximation in $L^1$ of trigonometric polynomials. By a trigonometric (correspondingly, analytic) polynomial of degree $n$ we mean a linear combination of harmonics $z^k$, $|k| \le n$ (correspondingly, $0 \le k \le n$). Every trigonometric polynomial can be regarded as a rational function with multiple pole at the origin. Hence, the result below can be readily obtained from Theorem \ref{mrt}. However, we would like to give a separate proof as an example of using Theorem \ref{t1}. 
\begin{Prop}\label{p4}
Let $g$ be a trigonometric polynomial of degree $n \ge 1$ and let $p_g$ be its best $H^1$--approximation. Then $p_g$ is an analytic polynomial of degree at most~$n$. Moreover, the function $g - p_g$ has the form 
\begin{equation}\label{e10}
g - p_g = const  \cdot \bar  z^n \prod_{1}^{K}(1 - \bar \lambda_k z)(z - \lambda_k)\prod_{1}^{M}(1 - \bar \mu_m z)^2,
\end{equation}
where $|\lambda_k|< 1$, $|\mu_m| \le 1$, and $K + M \le n-1$.
\end{Prop}
\beginpf Consider the inner function $\theta_n = z^n$. By the assumption,  $g \in \bar \theta_n H^1 \cap \theta_n \ov{H^1}$. The coinvariant subspace $K^{2}_{\theta_n}$ consists of analytic polynomials of degree at most $n-1$.  It follows from Theorem \ref{t1} that $g - p_g = \bar z^{n} I F^2$, where $F$ is an analytic polynomial of degree at most $n-1$ and without zeroes in the open unit disk; $I$ is a finite Blaschke product such that $IF$ is an analytic polynomial of degree at most $n-1$. Denote by $\lambda_k$ the zeroes of $I$ and by $1/\bar \mu_m$ those zeroes of $F$ that are not poles of $I$, taking into account multiplicities.
 It is now evident that the function $g - p_g$ is of form \eqref{e10}. Since $g$ and the right side in \eqref{e10} are trigonometric polynomials of degree at most $n$, the function $p_g$ is an analytic polynomial of degree at most~$n$.
\qed

\medskip 

\noindent {\bf Proof of Theorem \ref{mrt}.} Let $1 \le p \le \infty$, and let $g$ be a rational function with $n$ poles $\beta_i$ in the open unit disk, each counted according to multiplicity. Then   $g = h/B$, where $h \in H^p$ and $B$ is the Blaschke product with zeroes $\beta_i$, 
$$
B = \prod_{i=1}^{n}\frac{z - \beta_i}{1 - \bar\beta_i z}.
$$
On the unit circle $\T$ we have $g = \bar B h$. Let $p_g$ denote best $H^p$--approximation of $g$. By Theorem \ref{t1} (see also Remark 2.1 for the case $p=\infty$), the function $g - p_g$ can be uniquely represented in the form $g - p_g = c \bar B I F^{2/p}$, where $F$ is an outer function in $K_B^2$ and $I$ is an inner function such that $IF \in K_{B}^{2}$. 

It follows from the definition of $K^2_B$ that every function $f \in K_{B}^{2}$ has the form $P_f/Q_B$, where $Q_B=\prod_{i=1}^{n} (1 - \bar \beta_{i} z)$ and $P_f$ is an analytic polynomial of degree at most $n  - 1$. Since the function $F$ is outer, the polynomial $P_F$ has no zeroes in the open unit disk. Let us write it in the form $P_F = c_1 \cdot \prod_{1}^{n-1} (1 - \bar \alpha_{i} z)$, where $c_1$ is a constant and $|\alpha_i| \le 1$ (if $\deg P_F< n-1$, we let some of $\alpha_i$'s equal to zero). By the construction, $IF \in K^2_B$. Hence,  we have $I = \prod{}^{'} \frac{z - \alpha_i}{1 - \bar \alpha_i z}$, where the product $\prod{}^{'}$ is extended over all, some, or none of the $\alpha_i$ with $|\alpha_i|<1$. This yields formula \eqref{e1}. 
The theorem is proved.  \qed 

\medskip

\noindent {\bf Remark 3.1.} The dual extremal function $h_g$ to the function $g$ has the form
\begin{equation}\notag
h_g = c_2 \cdot z  \prod{}^{''} \frac{z - \alpha_i}{1 - \bar \alpha_i z} \prod_{1}^{n-1} (1 - \bar \alpha_i z)^{2/p'}  \prod_{1}^{n} (1 - \bar \beta_i z)^{-2/p'},  \\
\end{equation}
where $\prod{}^{''}$ is complementary to $\prod{}^{'}$ with respect to the $\alpha_i$ with $|\alpha_i|< 1$ and $c_2$ is a constant. Indeed, this follows  from Remark 2.2.

\section{Proof of Theorem \ref{t2}}\label{s4}
A bounded analytic function $\theta$ in the upper halfplane $\C_+$ of the complex plane $\C$ is called inner if $|\theta|=1$ almost everywhere on the real line $\R$ in the sense of angular boundary values. Coinvariant subspaces of the Hardy space $\H^p$ in $\C_+$ have the form $\mathcal K^{p}_{\theta} =\H^p\cap \theta \ov{\H^p}$.
Theorem \ref{t1} holds for functions $g$ in $\bar \theta \H^p$, as can be easily seen from its proof. We will deduce Theorem \ref{t2} from the following more general result.
\begin{Prop}\label{p1}
Let $\theta$ be an inner function in $\C_+$ and let $g \in \ov{\theta} \H^1 \cap \theta \ov{\H^{1}}$. Then we have $p_g \in \mathcal K^{1}_{\theta}$ for best $\H^1$--approximation $p_g$ of $g$. 
\end{Prop}
\beginpf By Theorem \ref{t1}, we have $g - p_g = \bar \theta I F^{2}$, where $F$, $IF$ are functions in $\mathcal K^{2}_{\theta} $. Hence, the function $g - p_g$ belongs to the subspace 
$$\bar \theta \cdot (\H^2 \cap \theta \ov{\H^2})\cdot (\H^2 \cap \theta \ov{\H^2}) \subset \bar \theta \cdot (\H^1 \cap \theta^2 \ov{\H^1}) \subset \bar \theta \H^1 \cap \theta \ov{\H^1}.$$
It follows that the function $p_g$ lies in the subspace $\mathcal K^{1}_{\theta} = \H^1 \cap \theta \ov{\H^1}$. \qed

\bigskip

\noindent {\bf Proof of Theorem \ref{t2}.} Consider the inner function $S^{a}: z \mapsto e^{iaz}$ in the upper halfplane~$\C_+$. A function $f$ in $L^1(\R)$ belongs to the Hardy space $\H^1$ if and only if $\supp \hat f \subset [0, +\infty)$. It follows that
every function $g \in L^1(\R)$ with $\supp \hat g \subset [-a, a]$ belongs to the subspace $\ov{S^a} \H^1 \cap S^a \ov{\H^{1}}$. By Proposition~\ref{p1}, we have $p_g \in \H^1 \cap S^a \ov{\H^{1}}$. Hence, $\supp \hat g \subset [0, a]$ and the result follows. \qed

\section{Interpolation problems related to best analytic approximation}\label{s5}
The problem of best $H^p$--approximation for functions in $\bar \theta H^p$ can be rewritten in the following form: given a function $g \in H^p$, find a function $h \in H^p$ such that the norm $\|g - \theta h\|_{L^p}$ is minimal. This is the problem of {\it constrained interpolation in $H^p$} with respect to the inner function $\theta$. 
An account of results related to constrained interpolation in $H^\infty$ is available in Chapter 3 of \cite{MR1864396}. We will consider the same problem in $H^p$, $1 \le p \le \infty$. Our observations are in line with \cite{MR0036314}, where the problem of best analytic approximation in $L^p$ for rational functions is reduced to a problem of interpolation.
For $\theta$ is an inner function, define the class $E_{\theta,p}$ by the formula
\begin{equation}\label{e7}
E_{\theta,p} = \{cIF^{2/p}: \; c \in \C, \;  I \mbox{ is inner, }\; F \mbox{ is outer,}\; IF \in \Kth\}.
\end{equation}
If $p$ is finite, there is no need in the constant $c$ in \eqref{e7}. We say that a function $f_2 \in H^p$ interpolates a function $f_1 \in H^p$ with respect to the inner function $\theta$ if $f_1 - f_2 \in \theta H^p$. For example, $f_1$ interpolates $f_2$ with respect to $z^n$ if and only if $f_{1}^{(k)}(0) = f_{2}^{(k)}(0)$ for all integers $0 \le k \le n-1$. Another example: for $\theta$ is a Blaschke product with simple zeroes $\Lambda$, $f_1$ interpolates $f_2$ with respect to $\theta$ if and only
if $f_1(\lambda) = f_2(\lambda)$ for all $\lambda \in \Lambda$. 

\medskip

The main result of this section is following.

\begin{Prop}\label{p2}
Let $1 \le p < +\infty$ and let $\theta$ be an inner function. 
Each function $f_1 \in H^p$ can be interpolated by a unique function $f_2 \in E_{\theta, p}$ with respect to $\theta$. Moreover, we have $\|f_2\|_{L^p} = \dist_{H^p} (f_1, \theta H^p)$.
\end{Prop}
\beginpf Take a function $f_1 \in H^p$ and set $g = \bar \theta f_1$. Let $p_g$ denote best $H^p$--approximation of $g$. By Theorem \ref{t1}, we have $g - p_g = c \cdot \bar \theta IF^{2/p}$, where the function $f_2 = c \cdot IF^{2/p}$ belongs to the class $E_{\theta, p}$. Note that $f_1 - f_2 \in \theta H^p$. Hence, the function $f_2$ interpolates $f_1$ with respect to the inner function $\theta$. By the construction, we have $\|f_2\|_{L^p} =   \dist_{H^p} (f_1, \theta H^p)$. 

Let us now prove that the interpolating function $f_2$ is unique. Suppose that there is an another function $f_2^* = c^*I^* {F^{*}}^{2/p}$ in $E_{\theta, p}$ that interpolates $f_1$ with respect to~$\theta$. We may assume that functions $F$, $F^*$ are of unit norm in $\Kth$ and have positive values at the origin. Let also $c >0$, $c^* > 0$. Consider the inner function $J$ such that $IJF = \bar z \theta \bar F$. Since $cIF^{2/p} - c^*I^*{F^{*}}^{2/p}$ lies in $\theta H^p$, we have
\begin{equation}\label{e5}
c = \int_\T c \bar \theta IF^{2/p} \cdot zJF^{2/p'} \, dm = \int_\T c^* \bar \theta I^*{F^{*}}^{2/p} \cdot zJF^{2/p'}\, dm \le c^*.
\end{equation}
Symmetric argument tells us that $c_* \le c$. Hence, we have equality in \eqref{e5}. It follows that outer functions $F$, $F^*$ have the same modulus on $\T$. Since $F(0)> 0$ and $F^*(0) > 0$, we have $F=F^*$. Again by equality in \eqref{e5}, inner functions $I$ and $I^*$ have the same argument on $\T$. Hence, $I=I^*$ and $f_2 = f_{2}^{*}$. \qed 

\medskip

A function $g \in L^p$ is called $H^p$--badly approximable if the zero function is the best analytic approximation of $g$ in $L^p$. Theorem \ref{t1} and Proposition \ref{p2} allow us to describe all $H^p$--badly approximable functions in $\bar \theta H^p$, where $\theta$ is an inner function and $1 \le p < \infty$.
\begin{Prop}\label{p3}
Let $1 \le p< \infty$. A function $g \in \bar \theta H^p$ is $H^p$--badly approximable if and only if $\theta g \in E_{\theta, p}$.
\end{Prop}
\beginpf By Theorem \ref{t1}, we have $\theta g \in E_{\theta, p}$ for every $H^p$--badly approximable function $g \in \bar \theta H^p$. Conversely, take a function  $g \in \bar \theta E_{\theta, p}$ and consider its best $H^p$--approximation $p_g$. Set $f_1=\theta g$. The function $f_2 = f_1 - \theta p_g$ interpolates $f_1$ with respect to the inner function $\theta$. By Theorem \ref{t1}, we have $f_2 \in E_{\theta, p}$. Hence, two functions $f_1, f_2 \in E_{\theta, p}$ interpolate the function $f_1$ with respect to $\theta$. It follows from Proposition \ref{p2} that $f_1 = f_2$ and so $p_g =0$. \qed 

\medskip

We conclude this section with some examples.

\medskip

\noindent {\bf Example 1.} The class $E_{z^n, 1}$ consists of polynomials of the form
\begin{equation}\label{e8}
const \cdot \prod_{1}^{K}(1 - \bar \lambda_k z)(z - \lambda_k)\prod_{1}^{M}(1 - \bar \mu_m z)^2,
\end{equation}
where $|\lambda_k|< 1$, $|\mu_m| \le 1$, and $K + M \le n-1$. As to the author knowledge, the problem of {\it constructive} interpolation by polynomials of form \eqref{e8} with respect to~$z^n$ is open until now. A detailed discussion of this problem can be found in Section 5 of~\cite{MR0036314}.

\medskip

\noindent {\bf Example 2.} Let us compute the upper bound of quantities $|f(0) + f'(0)|$ over all $f \in H^1$ of unit norm. By duality and Proposition \ref{p2}, this problem reduces to interpolation of $1 + z$ with respect to $z^2$ by a polynomial in $E_{z^2, 1}$. It is easy to see that the polynomial $\frac{1}{4}(2+z)^2$ solves this problem. Hence, we have 
$$\sup\left\{|f(0) + f'(0)|: \; f \in H^1, \|f\|_{H^1} = 1\right\} = \frac{1}{4}\int_\T|z+2|^2\, dm = \frac{5}{4}.$$ 
The general problem of calculation of $\sup\left\{|a_0 f(0) + a_1 f'(0)|: \; f \in H^1, \|f\|_{H^1} = 1\right\}$ is solved in Section 5 of~\cite{MR0036314}.
 
\medskip

\noindent {\bf Example 3.} The problem of constrained interpolation in $H^\infty$ with respect to the inner function $z^n$ is the classical Schur problem. It can be stated as follows: given $n$ numbers $a_0, \ldots a_{n-1}$, find a function $f \in H^\infty$ of minimal norm such that $f^{(k)}(0)/k! = a_k$ for all $k$. This problem can be solved constructive, see \cite{Schur17} or Paragraph 3.4.2.(ii).(b) in \cite{MR1864396}. The solution has the form 
\begin{equation}\label{e9}
const \cdot \prod_{k=1}^{K} \frac{z - \lambda_k}{1 - \bar \lambda_k z}, \qquad |\lambda_k| < 1, \quad K \le n-1. 
\end{equation}
Clearly, class $E_{z^n, \infty}$ consists of functions of form \eqref{e9}. This agrees well with Proposition \ref{p2} (its proof works in the case $p = \infty$ if there exists the dual extremal function to the function $\bar \theta f_1$). Similarly, the classical Nevanlinna-Pick problem reduces to interpolation by functions of form \eqref{e9} with respect to a finite Blaschke product. For more information, see Chapter 3 in \cite{MR1864396}. 
\bibliography{bibfile}
\bibliographystyle{plain} 
\enddocument